\documentclass[11pt]{amsart}
\usepackage{amsmath,amssymb, graphics, amscd,latexsym}
\makeatletter

\newtheorem{Theorem}{Theorem}
\newtheorem{Lemma}[Theorem]{Lemma}

\newtheorem{Proposition}[Theorem]{Proposition}

\newtheorem{Definition}[Theorem]{Definition}
\newtheorem{Remark}[Theorem]{Remark}

\newcommand\al{\alpha}

\newcommand\si{\sigma}

\newcommand\Si{\Sigma}

\newcommand\cM{\mathcal  M}

\newcommand\cI{\mathcal  I}

\newcommand\cO{\mathcal  O}

\newcommand\bfC{\mbox {\bf  C}}

\newcommand\bfP{\mbox {\bf  P}}

\newcommand\bfQ{\mbox {\bf  Q}}

\newcommand\bfZ{\mbox {\bf  Z}}

\newcommand\order{\mbox{\text{order}\/}}

\newcommand\PGL{\text{PGL}\/}

\newcommand\ord{\text{ord}\/}

\newcommand\Coeff{\text{Coeff}\/}

\begin{document}
\title[Zariski pairs on sextics I 
]
{  Zariski pairs on sextics I }

\author
[M. Oka ]
{Mutsuo Oka }
\address{\vtop{
\hbox{Department of Mathematics}
\hbox{Tokyo  University of Science}
\hbox{1-3 Kagurazaka,Shinjuku-ku}
\hbox{Tokyo 162-8601}
\hbox{\rm{E-mail}: {\rm oka@ma.kagu.tus.ac.jp}}
}}
\begin{abstract}
We study Zariski pairs of sextics which are distinguished by
the Alexander polynomials. For this purpose,
we present two constructive methods to produce explicit sextics of
 non-torus type
with  given configuration of simple singularities. 
\end{abstract}

\maketitle

\pagestyle{headings}
\section{Introduction}

A plane curve of degree $d$ is called a  $(p,q)$-torus curve ($p,q\ge 2$)
or a curve of $(p,q)$-torus type if $p,q $  divide $d$ and 
$C$ can be defined by a polynomial of the form
$f_{d/p}(x,y)^p+f_{d/q}(x,y)^q$ where the degree of $f_{d/p},\, f_{d/q}$
are $d/p$ and $d/q$ respectively. 
In this paper, we consider  sextics
of (2,3)-torus type. Their defining polynomial $f(x,y)$ are written as
$f=f_2^3+f_3^2$.  A curve of  $(p,q)$-torus type has 
$d^2/(pq)$ cusp singularities of the type $y^p+x^q=0$
at the intersection
$f_{d/p}=f_{d/q}=0$ if  two curves $f_{d/p}=0$ and $f_{d/q}=0$ intersect
transversely. Then the global Alexander polynomial of $C$ is the same with the
local
Alexander polynomial of the torus  link in  $S^3$ defined by 
$y^p+x^q=0$ ( \cite{OkaSurvey}).
A curve $C$ of degree $d$ is called {\em of non-torus type} if the defining 
polynomial has no such expression for any $p,q\ge 2$ with $p|d$ and
$q|d$. 
For a given torus curve $C$ and 
 a singularity $P\in C$, we say that $P$ is
an {\em inner} singularity if it is on the intersection $f_{d/p}=f_{d/q}=0$.
A torus curve $C$ is called {\em tame} 
if it has only inner singularities. 

Two plane curves $C,C'$ of the same degree $d$ is called a {\em  Zariski pair}
if  (i) they have the same configuration of the singularities and there exist
topologically equivalent tubular neighborhood, and
(ii) the topology of the pair $(\bfP^2,C)$ and $(\bfP^2,C')$ are
not  equivalent. 

This notion is introduced by Zariski  \cite{Za1} and formulated as above by
Artal Bartolo \cite{Artal}.
The first such example  is given by Zariski
on sextics:  $(Z_6,Z_6')$ where two curves $Z_6,Z_6'$
are sextics with 6 $A_2$-cusps,  $Z_6$ is a sextic of (2,3)-torus type
but $Z_6'$ is not of torus type \cite{Za1}. The corresponding Alexander polynomials
are
$t^2-t+1$ and $1$ respectively.
 The purpose of this paper is to give further Zariski pairs in sextics
 which can be distinguished by Alexander polynomials.
D.T. Pho  determined all possible configurations of the singularities
of irreducible tame sextics of (2,3)-torus type 
$C: \, f_2^3+f_3^2=0$ with simple singularities (\cite {Pho}).
Hereafter, we simply say a ``configuration'' instead of a ``configuration
of the  singularities''.
In fact, they are given by
\begin{enumerate}
\item (1,1,1,1,1,1): \,$[6A2]$\item  (2,1,1,1,1):\, $[A_5,4A_2],\,[E_6,4A_2]$
\item (2,2,1,1): $[2A_5,2A_2],\,[A_5,E_6,2A_2],\, [2E_6,2A_2]$
\item (2,2,2):\,$ [3A_5],\,[2A_5,E_6],\,[A_5,2E_6],\,[3E_6]$
\item (3,1,1,1):\, $[A_8,3A_2]$\item
      (3,2,1):\,$[A_8,A_5,A_2],\,[A_8,E_6,A_2]$\item (3,3):\,$[2A_8]$
\item (4,1,1):\,$[A_{11},2A_2]$\item (4,2):\,$ [A_{11}, A_5],\,[A_{11},E_6]$
\item (5,1): \,$[A_{14},A_2]$\item  (6):\, $[A_{17}]$
\end{enumerate}
The numbering corresponds to the partition types of $6$ (=the total intersection number
of the conic $f_2=0$ and the cubic $f_3=0$).
In this paper, we will show that 
\begin{Theorem}
For any configuration $\Si$ in (1)$\sim$ (11), there is a Zariski partner
irreducible sextic, say $C'$,
of non-torus type with configuration 
$\Si$. So choosing an irreducible sextic of torus type $C$
with  configuration $\Si$, we have a Zariski pair
$(C,C')$ and they are distinguished by the 
Alexander polynomial, either $t^2-t+1$ or $1$ respectively.
\end{Theorem}
Most of them has been studied by several poeple.
We have studied the fundamental group of the complement of the sextics
with configuration
$[6A_2],\, [3E_6]$ in \cite{OkaSymmetric,Okadual}. 
Tokunaga studied the cases with $E_6$ singularity in \cite{Tokunaga-Some1,Tokunaga-Some2}.
 In \cite{Eyral-Oka1}, Eyral and the author have 
proved the fundamental groups
of the complements of the  sextics of non-torus
 type with  the
following configurations are abelian.
\[
 [2A_8],\,[A_{17}],\,[A_{11},E_6],\,[A_{11},A_5],\,[A_{14},A_2],\,[A_{8},A_5,A_2],\,[A_8,E_6,A_2]
\]
Sextics with  $A_{17}$ is discussed in \cite{Artal}. See also \cite{INO}.
So new results in this paper are  only those assertions
 for the remaining cases.
Especially for the cases $[A_5,4A_2], \, [E_6,4A_2]$, the construction
is a new result. We also construct 4 dimensional subspace of
sextics of  non-torus type with $6A_2$ and the fundamental group of the 
complements are abelian.

 In the second part of this paper, we will study Zariski
pairs in non-irreducible sextics \cite{ZP2}.

 \section{Two computational methods}
 In this section, we introduce two computational methods to 
 obtain a sextics with a prescribed configuration of singularities.
 First we recall basic properties of sextics of torus type.
 Consider a sextic  of torus type
 \[T:\, f(x,y)=f_2(x,y)^3+f_3(x,y)^2=0\]
 where $f_2, \,f_3$ are polynomials of degree $2$ and $3$ respectively.
 Hereafter, we  mean by a sextic  of torus type a sextic of (2,3)-torus type.
 Recall that a  singularity $P\in T$
 is called an inner singularity if 
 $P$ is at the intersection of the conic $C_2:\,f_2(x,y)=0$
 and the cubic $C_3:\, f_3(x,y)=0$. If $C_3$ is non-singular at 
 $P$, the corresponding singularity of $T$ at $P$ is  described
 by the local intersection number $\iota:=I(C_2,C_3;P)$ and 
 the  singularity at $P$ is isomorphic to 
 $A_{3\iota-1}$. If $C_3$ is singular at $P$, $(T,P)$ is a simple
 singularity if and  only if when $C_2$
 is reduced (so smooth at $P$) and $I(C_2,C_3;P)=2$.
 In this case, we have $(T,P)\cong E_6$.

 In \cite{OkaAtlas}, we have  introduced an local invariant
  $\rho(P,5)$ for any local singular point $P\in C$
 by $\dim_{\bf C}\cO_P/\cI_{P,5,6}$ where $C$ is a sextic and the ideal
 $\cI_{P,5,6}$ is defined by the
 adjuction ideal (\cite{Libgober3,Loeser-Vaquie,Esnault, Artal}) and put
 $\rho(C,5):=\sum_{P}\rho(P,5)$.
 For a sextic of torus type, we observe (\cite{OkaAtlas}).

 \begin{Proposition} Assume that $T$ is a sextic of torus type and 
 $P$ be an inner simple singularity. Then
 $\rho(P,5)=\iota(C_3,C_2;P)$.
 Thus 
  $\sum_{P:inner} \rho(P,5)=6$ for a sextics of torus type
 with only simple singularities where the summation is taken for all inner
 singularities $P$.
 \end{Proposition}


 \vspace{.3cm}
 Suppose that  a sextic $C=\{g(x,y)=0\}$
 is given and assume that the singularities of $C$ are simple singularities
  contained in the list:
$ \{A_1,\,E_6,\,A_2,\,A_5,\,A_8,\,A_{11},\, A_{14},\, A_{17}\}$.
 Here $C$ need not to be of torus
 type.
 Let us denote the configuration set of singularities of $C$
  by
 $\Si(C)$.
  The following is 
 a convenient criterion by H. Tokunaga to check if $C$ is of torus type.
 \begin{Lemma}{\bf (Tokunaga Criterion \cite{Tokunaga-torus})}\label{Tokunaga-Criterion}
 The sextic $C$ is of torus type if and only if   there exists a conic
 $C_2: g_2(x,y)=0$
 such  that $C_2\cap C\subset \Sigma(C)$ and 
 $I(C,C_2;P)=2\rho(P,5)$ for any  $P\in C\cap C_2$.
 \end{Lemma}

 This is a useful tool to check if  a given sextic $C$ is of torus
 type or not when $C$ is already given. For our purpose,
 we are interested in explicit construction of  sextics of 
 non-torus type
 with given configuration of singularities $\Si$.

 \begin{Definition}
 Assume that $O\in C$ is an  $A_{3\iota-1}$ singularity.
 We assume that $y=t_1x$ is the tangent cone of $(C,O)$.
 Then there exist unique complex numbers 
  $t_2,\dots,t_\tau$, with $\tau=[3\iota/2]$
 such that we can write
 \[\tilde g (x,y):=g(x,y+t_1x+t_2x^2+\dots +t_\tau x^\tau)
 =a y^2+ b x^{3\iota} +\text{(higher terms)},\, a,b\ne 0
 \]
 Higher terms are linear combinations of monomials $x^\alpha\,y^\beta,\, 2\alpha+3\iota \beta>6\iota$.
 We call $y=t_1x+t_2x^2+\dots+t_\tau x^\tau$ {\bf the maximal contact
 coordinate curve} of $(C,O)$.
 \end{Definition} Let $E$ be the curve defined by 
 $y-(t_1x+t_2x^2+\dots+t_\tau x^\tau)=0$. Then by definition, we have
 $I(C,E;O)=\order_x\, g(x,t_2x^2+\dots+t_\tau x^\tau)=3\iota$.

 \vspace{.5cm}
 {\bf Method 1.  $\sharp(\Si)\le 4$.}
 By the action of $\PGL(3,\bfC)$, we may assume that the locations 
 of singularities are  on the chosen explicit points.
 For example, we explain the construction of a sextic of non-torus type
 with  $\Si=[A_{14}+A_2+A_1]$.
 We start from the generic form of sextic:
 \[
  g(x,y)=\sum_{i+j\le 6} a_{ij} x^i\,y^j
 \]
 We have 27 free coefficients.
 Let us assume that $C=\{g(x,y)=0\}$
 has $A_{14}$ singularity at $O=(0,0)$, $A_2$ at $P=(1,0)$
 and $A_1$ at $Q=(-1,1)$. The condition $(C,Q)\cong A_1$ is given by 3
 linear equation. The condition $(C,P)\cong A_2$ is given by 3 linear
 condition and one quadratic equation i.e., the hessian vanishes at $P$.
 To avoid this quadratic equation,
 we may assume that $x=1$ is the tangent cone at 
 $P$ and  $x=0$ is the tangent cone at the origin. Thus the quadratic
 equation is replaced by two linear equation
 $\frac{\partial^2 g}{\partial x \partial y}(P)=\frac{\partial^2
 g}{\partial y^2}(P)=0$.
 So up to here, we eliminated 8 coefficients.
 The most difficult one is the condition for $O$ to be $A_{14}$.
 For this purpose, we consider the maximal contact coordinate curve for $A_{14}$.
 Namely we suppose that $x=0$ is the tangent cone of $(C,O)$ and 
   $x=t_2\, y^2+\dots+t_7 \, y^7$ is
 the maximal contact coordinate curve at $O$. Here $t_2,\dots, t_7$ are
 variables to be determined later.
 Thus  to find a sextic of non-torus type $C$,
 we consider the polynomial
 $h(x,y):=g(x+t_2\, y^2+\dots+t_7 \, y^7,\, y)$
 and put it as $h(x,y)=\sum _{ij} b_{ij} x^i \, y^j$.
 (This has 19+6=25  free variables in their coefficients, where 19 of
 them are linear.)
 We use the notation ${\rm Coeff}\,(h, x^i \, y^j)$ for $b_{ij}$.
 Note that $b_{ij}$ is a linear combination of $a_{k\ell}$ with 
 coefficients in $\bfC[t_2,\dots, t_7]$.
 The assumption \[
 h(x,y)=a\,x^2+c\,y^{15}+\text{(higher terms)},\quad a,c\ne 0
		\]
 implies that
 \[
   b_{0j}=0,\, j\le 14\quad
 \text{and}\,\,\, b_{1s}=0,\, s\le 7.
 \]
 which is comparatively easy to be solved.
 We have 23 equations. So at the end of calculation, we have to solve
 some non-linear equations. 
 Now we explain how we can solve the equations explicitly
 (preferably over $\bfQ$)
 so that the resulting sextics are of non-torus type.
 Assume that $C$ is a sextic of torus type i.e., there exists a conic 
 $C_2:=\{h_2(x,y)=0\}$ and a cubic
 $C_3:=\{h_3(x,y)=0\}$ such that $g(x,y)=h_2(x,y)^3+h_3(x,y)^2$. Then 
 $C_3$ is smooth at $O$ and assume that the equation $h_3(x,y)=0$ can be solved
 locally as 
 $x=\sum_{i\ge 2} s_2\, y^i$.  Put $x_1=x-\sum_{i\ge 2} s_2\, y^i$.
 We see that the assumption  $I(C_2,C_3;O)=5$ 
 implies that $\ord_y\,h_2(\sum_{i\ge 2}s_i \, y^i,y)=5$.
 This implies that $h_3=x_1\, U$ with $U$ being a unit,
 $h_2(x_1+\sum_{i\ge 2}s_i \, y^i,y)= \al x_1+ \beta\, y^5+\text{(higher
 terms)},\,\al,\,\beta\ne 0$ and 
 \[
  h_2(x_1+\sum_{i\ge 2}s_i \, y^i,y)^3+h_3(x_1+\sum_{i\ge 2}s_i \, y^i,y)^2=
 U(0) x_1^2+\beta^3 \, y^{15}+\text{(higher terms)}
 \]
 Thus by the uniqueness of the maximal contact coordinate curve, we have
  $ t_i=s_i,\,i\le 5$. Thus
 {\em The maximal contact coordinate of inner $A_{3\iota-1}$
 singularity of a sextic of torus type is given by solving $f_3(x,y)=0$
 in $x$.}
%
 Solving the five equalities
 \[
  h_2(1,0)=0,\quad \text{and}\,\,
  \Coeff(h_2(t_2y^2+\dots +t_5y^5,\,y),\,y^j)=0,\quad j=0,\dots, 4
 \]
 and thus  we see that $h_2(x,y)$ is written as
 $h_2(x,y)=c\, (t_2^2y^2+t_3xy+x^2t_2-xt_2),\, c\ne 0$
 and  the last equality
 $  \Coeff(h_2(t_2y^2+\dots +t_5y^5,y),y^5)=0$ implies
 the equality
 \[
  J:=t_3t_4+2t_2^2t_3-t_5t_2=0
 \]
 Note that there are 5 free coefficient in $h_2$ and 6 equations can not 
 satisfied without the condition $J=0$.
 This implies that the maximal contact coordinate for a sextic of a torus
 type
 must satisfy $J(t_2,\dots, t_5)=0$,
 under the  above assumption on the location of the singularities.
  Thus in the construction of sextics
 with configuration $[A_{14}+A_2+A_1]$,
 we choose the coefficients $t_2,\dots,t_5$ such that $J(t_2,\dots, t_5)\ne 0$
 which is enough to guarantee that the sextic to be obtained will be of non-torus type.
 An example  of sextic with $[A_{14}+A_2+A_1]$  is give by:
 \begin{multline}
 g(x,y)=
  208716480\,x^5-74370240\,x^6-38384640\,y^2\,x+331067520\,y^4\,x+19192320\,x^2\\
 -660935520\,y^2\,x^4
 +19192320\,y^4-175129920\,x^4-701719200\,y^4\,x^2-9596160\,y\,x^2+28788480\,y^5\\
 +760495680\,y^3\,x^2-226709280\,y^2\,x^2
 -106757280\,y^6-319072320\,y\,x^3-374250240\,y^5\,x\\
 +926029440\,y^2\,x^3-338264640\,y\,x^5+666933120\,y\,x^4\\
 -19192320\,y^3\,x-858856320\,x^3\,y^3+21591360\,x^3
 \end{multline}

 {\bf Method 2: Semi-torus  method}.  This method is quite useful for the
 case
 $\sharp (\Si)\ge 5$ but it is also convenient for the case $\sharp
 (\Si)\le 4$.
 In the case $\sharp (\Si)\ge 5$, the method  1 does not work in general as 
 the locations of all singularities can not be arbitrary, 
 though we can   fix the locations of four singularities.
 This is the case when we consider the configurations
  for example $4A_2+A_5,\, 4A_2+E_6,\, 6A_2$. 
 To consider these cases systematically, 
 we consider {\em a sextic  of semi-torus 
 type} which is defined by
 \[
 (ST)\quad C:\, g(x,y)\,=\,f_2(x,y)^3+g_2(x,y)^2\,h_2(x,y)=0
 \]
 where $f_2,g_2$ (respectively $h_2$) are polynomials of degree 2
 (resp. of degree $\le $ 2).
 A singularity $P\in C$ is called {\em wild} if $P$    is on the
 intersection of three conics:
 $ \{f_2(x,y)=g_2(x,y)=h_2(x,y)=0\}$.
 We say $C$ is {\em nice} if $C$ does not have any wild singularities.
 As in the case of (2,3)-torus curve, we call a  singular point
 $P\in C$  {\em inner} (respectively {\em outer}) if 
 $f_2(P)=g_2(P)=0$ (resp. $f_2(P)\ne 0$) and $P$ is not wild.
 By the same argument as in  \cite{Pho}, 
 it is easy to see that 
 the possible configurations of inner singularities
 of a nice sextic $C$ of semi-torus type are 
 \[
 (\sharp):\,
  [4A_2],\, [A_5+2A_2],\, [E_6+2A_2],\, [2E_6],\, [A_8,A_2],\, [A_{11}]
 \]
 We can check by a direct computation
 that {\em  the  possible outer singularities are $A_j,\, j\le 5$ and
 $E_6$ is not possible as an outer singularity}. (We do not consider other singularities.)
 So we can construct sextics with configuration
 $\Si=\Si_{inner}\cup\Si_{outer}$ where $\Si_{outer}$ can be 
 either $[2A_2]$ or $ [A_5]$.
 Just like the case of sextic of torus type, we can show easily that
 \begin{Proposition}
 An inner singularity $P\in \{f_2(x,y)=g_2(x,y)=0,\,h_2(x,y)\ne 0)\}$ is 
 an $E_6$ if  only if the conic $g_2(x,y)=0$ is singular at $P$ and
  $I(f_2,g_2;P)=2$. This implies that 
 $g_2(x,y)=0$  is either a union of two lines meeting at $P$ or 
 a line of mutiplicity $2$ passing through $P$.
 \end{Proposition}
  Thus  $[2E_6]$ is obtained only if
 $g_2(x,y)=0$ is a double line $\ell(x,y)^2=0$ and 
 $\ell(x,y)=0$ intersects the conic $f_2(x,y)=0$ at two 
 distinct points.
 We assume that a sextic $C$ of semi-torus type is given in the form
 $(ST)$
 and its inner configuration $\Si_{inner}$ is given
  from the list $(\sharp)$.
 To explain a little more in detail, we   construct  sextics of semi-torus type
 with $\Si_{inner}\cup \{2A_2\},\, \Si_{inner}\cup \{A_5\}$
 systematically.
 Let $P,Q$ be  two outer $A_2$ singularities.
 In the case $P=Q$,  we assume that  $(C,P)\cong A_5$.
 First we put 
 $ t_1:=f_2(P),\, s_1:=g_2(P),\, t_2:=f(Q),\, s_2:=g_2(Q)$.
 We assume that $t_1,\,t_2,\,s_1,\,s_2 \ne 0$ and 
 \begin{eqnarray}
 &\,\, (C,P)\cong A_2,\,(C,Q)\cong A_2,\quad \text{if}\,  P\ne Q\\
 &\,\, (C,P)\cong A_5,\quad \text{if}\,\, P=Q
 \end{eqnarray}
 Then  we first observe
 that :
 \begin{Proposition} The linear system of conics $j_2=0$ satisfying
  $I(j_2,g;S)=2I(f_2,g_2;S)$
 for
 any inner singularity $S$ is given by
 \[
 C_2(t,s):\quad j_2(x,y):=t\, f_2(x,y)+ s\, g_2(x,y)=0,\quad s,t\in \bfC
 \]
 \end{Proposition}
 {\em Proof.} Note that the space of conics is 5-dimensional and 
 the condition  is 4-dimensional. So the space is a linear system.
  Assume that $S$ is an inner singular point. For a generic $t,s$,  the conic $j_2=0$ is
 smooth and 
 $I(j_2,g_2;S)= I(j_2,f_2;S)$, and therefore   we have
 $I(j_2,f_2^3+g_2^2h_2;S)=2I(j_2,g_2;S)$.\qed

 Now we assert
 \begin{Theorem} Assume that $C$ is a nice sextic  and let $\Si_{inner}$ be the 
 inner configuration of $C$. Then $C$ is a sextics of non-torus type 
 which contains

 {\rm (a)}
 $\Si_{inner}\cup\{2 A_2\}$, if $P\ne Q$ and $t_1,s_1,t_2,s_2\ne 0$
 and  $t_1s_2-t_2s_1\ne 0$, or

 {\rm  (b)} $\Si_{inner}\cup\{ A_5\}$, if $P=Q$ and $t_1,s_1\ne 0$
  and
 the tangent cone of $C$ at $P$
 is not the tangent line of the conic $C_2(s_1,-t_1)$ at $P$.
 \end{Theorem}
 {\em Proof.}
 Suppose that there exists a conic $j_2$ passing through  the singular points
  of $C$ with the prescribed intersection number. Then as it passes
  through $P,Q$, we
 must have the equality:
 $t_j\,t+s_j\, s=0$ for $j=1,2$.
  This has a solution $(t,s)\ne (0,0)$ if and only if $t_1s_2-t_2s_1=0$,
 contradicting to the assumption.
 If $P=Q$, the condition that $j_2=0$ passes through
 $P$ implies that $(t,s)\sim (s_1,-t_1)$. If $(C,P)\cong A_5$
 and if $C$ is of torus type  defines as  $j_2^3+j_3^2=0$ for some cubic
 form $j_3$, the
 tangent cone is given by 
 tangent line of $j_2=0$ at $P$.\qed
 \begin{Remark}
 1. The equation $f(P)=0$ is usually difficult to solve (i.e., the
  elimination
of a parameter involves a big denominator). This  difficulty is
  resolved by
 introducing new parameter $t_1,s_1$ so that $f_2(P)=t_1$ and $g_2(P)=s_1$
 as the latter
  are easy to be solved in a coefficient of $f_2$ or $g_2$ respectively. Then
  the equation
 $f(P)=0$ is a linear equation in the coefficients of $h_2$.

 \noindent
 2. For a wild singularity $P$ of a semi-torus sextic  $C$, so
 $f_2(P)=g_2(P)=h_2(P)=0$ with $\iota=I(f_2,g_2;P)$,
 the corresponding singularity is 
 $D_{3\iota+1}$ for $\iota=1,2,3,4$.
 \end{Remark}
 Let $\cM_{ST}(6A_2;6)$ be the configuration space of sextics with semi-torus type
 with $6A_2$ and the non-torus condition
 $t_1s_2-t_2s_1\ne 0$.
  This is a subspace of the configuration space of irreducible sextics with
 non-conical $6A_2$ (let us denote this space by
 $\cM_{nt}(6A_2;6)$  and we claim that
 {\em $\cM_{ST}(6A_2;6)/\PGL(3,\bfC)$ is connected of dimension $4$}.
 To see the irreducibility, we start from the generic forms of $f_2,\,g_2,\,h_2$
 and consider the slice condition: 

  $(\sharp)$:~~~2 outer $A_2$ are at
  $P_1:=(0,1),\,P_2:=(0,-1)$
  with the respective tangent cones $y\mp 1=0$.

  The
 corresponding conditions can be written as 
 \begin{eqnarray*}
 f_2(P_i)=t_i,\, g_2(P_i)=s_i,\,f(P_i)=0,\,i=1,2\\
 \frac{\partial f}{\partial x}(P_i)=\frac{\partial f}{\partial
  y}(P_i)=\frac{\partial^2 f}{\partial x\partial y}(P_i)
 =\frac{\partial^2 f}{(\partial x)^2}(P_i)=0,\, i=1,2
 \end{eqnarray*}
 Using these equations, we can eliminate 10 coefficients and the number of the 
 remaining  {\em free
 coefficients} are $17-10-1=6$. (One dimension comes from 
  the change 
 $(g_2,h_2)\mapsto (\alpha\,g_2,h_2/\alpha)$). This slice is
 clearly connected (or  irreducible) as  it is a Zariski open subset of 
 $\bfC^6$. As we have two dimensional isotopy group
 which fix $P_1,P_2$ and the tangent cone $y\mp 1=0$. So the quotient is
 4 dimensional. 
 \begin{Proposition}
 For any sextic $C\in \cM_{ST}(6A_2;6)$, the fundamental group 
 $\pi_1(\bfP^2-C)$ is isomorphic to the cyclic group
 $\bfZ/6\bfZ$.
 \end{Proposition}
   In fact by the connectivitity of the configuration slice,
  it is enough to show the assertion for one explicit 
   $C\in \cM_{ST}(6A_2;6)$.  Recall that  we have shown
   shown that $\pi_1(\bfP^2-C)\cong \bfZ/6\bfZ$ in \cite{OkaSymmetric} where
$C$ is defined by $f_6(x,y)=0$,
 \begin{multline*}
 {f_6(x,y)} := x^{2}\,(x - 1)^{2}\,(x^{2} + 2\,x) - x^{2}\,(x - 1
 )^{2}\,(y^{2} - 1) \\
 + {\displaystyle \frac {1}{3}} \,(x^{2} - 1)\,
 (y^{2} - 1)^{2} - {\displaystyle \frac {1}{27}} \,(y^{2} - 1)^{3}
 \end{multline*}
 On the other hand,  $f_6$ has a semi-torus decomposition as

 \begin{eqnarray*}
 &f_2(x,y)=- {\displaystyle \frac {1}{3}} \,y^{2} + x^{2} - {\displaystyle 
 \frac {2}{3}},\quad g_2(x,y)= - (x - 1)\,y,
 \quad
 h_2(x,y)=2\,x + 1
 \end{eqnarray*}
 \begin{Remark}
 We have a conjecture that the  space $\cM_{nt}(6A_2;6)$
 is connected and  the fundamental group of the complement
 $\pi_1(\bfP^2-C)$ is abelian for any $C$ in
  $\cM_{nt}(6A_2;6)$.
 The above Proposition  is a partial answer to this. We know  that the
  Alexander polynomial is trivial for any $C$ in
  $\cM_{nt}(6A_2;6)$.
 We expect that the dimension of 
  $\cM_{nt}(6A_2;6)/\PGL(3,\bfC)$
 is $19-6\times 2=7$ dimensional. We can check directly by computation
 that the dimension of the quotient space $\cM(\Si;6)/\PGL(3,\bfC)$ with a given configuration of simple 
 singularities $\Si$ is 19 minus the total Milnor number,
  if the number of singularities are less than 4.
 \end{Remark}
\section{Construction of irreducible sextics of non-torus type}
In this section, we construct irreducible non-torus sextics with configurations
(1), (2), (3), (4), (5) and (8).

\noindent
1. First we consider the configuration $\Si=[A_{11},2A_2]$.
This can be constructed either by Method  1 or also Method 2.
Usually the computation by Method 2 is less heavy.
For example, we look for a sextics of semi-torus type:
$$C:\,g(x,y)=f_2(x,y)^3+g_2(x,y)^2\,h_2(x,y)=0$$

where $A_{11}$ is at (1,0) as an inner singularity 
and we put
two $2A_2$ as  outer singularities at 
$(0,1),\,(0,-1)$. Foe example, we can take the equations:
\begin{multline*}
f_2 := 16962 + 43740\,y - 33924\,x - 9816\,\sqrt{3} - 
9816\,x^{2}\,\sqrt{3} - 25286\,y\,\sqrt{3} - 15792\,y^{2}\,\sqrt{
3} \\
\mbox{} + 19632\,x\,\sqrt{3} + 27445\,y^{2} + 16962\,x^{2} 
\end{multline*}
\[
g_2 := 6 + 6\,x^{2} - 38\,y^{2} + 9\,\sqrt{3} + 9\,x^{2}
\,\sqrt{3} - 11\,y^{2}\,\sqrt{3} + 16\,y\,\sqrt{3} - 18\,x\,
\sqrt{3} - 12\,x + 3\,y
\]
\begin{multline*}
h_2 :=  - ( - 136995863187\,x^{2} - 45018774416\,y^{2} + 
24984902205\,\sqrt{3} + 80265133314\,x^{2}\,\sqrt{3} \\
\mbox{} + 23477770096\,y^{2}\,\sqrt{3} + 49893109990\,y\,\sqrt{3}
 - 23083475220\,x\,\sqrt{3} \\
\mbox{} + 36411797484\,x\,y + 35951630964\,x - 83905166058\,y - 
41386138293 \\
\mbox{} - 20817696276\,x\,y\,\sqrt{3}) \left/ {\vrule 
height0.41em width0em depth0.41em} \right. \!  \! (( - 3539 + 
1787\,\sqrt{3})\,(7 + 4\,\sqrt{3})^{4})
\end{multline*}
The conics $f_2(x,y)=0$ and $g_2(x,y)=0$ are tangent at $(1,0)$ with
 intersection multiplicity 4 so that $C$ has $A_{11}$ singularity at
 $(1,0)$.

  \vspace{.5cm}\noindent
  2. Now we construct a sextics of non-torus type with
     $\Si=[A_8+3A_2]$
  as its configuration. We can use both methods in Case 1 or 2. Here we
     use the second method. So we consider sextic of semi-torus type
whose inner singularities are $A_8$ at (1,1) and two outer $A_2$ at
$(1,0),\,(0,1)$.
As an example, we obtain:
 \begin{eqnarray*}
 &g(x,y)=f_2(x,y)^3\,+\, g_2(x,y)^2\,h_2(x,y)\\
 &f_2 :=  - y^{2} - {\displaystyle \frac {1}{5}} \,y\,x\,
 \sqrt{59} - {\displaystyle \frac {13}{5}} \,y\,x + 
 {\displaystyle \frac {1}{5}} \,y\,\sqrt{59} + {\displaystyle 
 \frac {13}{5}} \,y  + 3\,x^{2} + {\displaystyle \frac {1}{5}} \,x
 \,\sqrt{59} - {\displaystyle \frac {12}{5}} \,x - {\displaystyle 
 \frac {1}{5}} \,\sqrt{59} + {\displaystyle \frac {2}{5}} \\
 &g_2 :=  - {\displaystyle \frac {1}{219010}} ( - 155 + 6\,
 \sqrt{59})(1810\,y^{2} + 173\,\sqrt{59} + 1724 + 1724\,y\,x - 
 5084\,y + 1506\,x \\
 &\mbox{} + 173\,y\,x\,\sqrt{59} - 233\,y\,\sqrt{59} - 1680\,x^{2}
  - 30\,x^{2}\,\sqrt{59} - 83\,x\,\sqrt{59}) \\
 &h_2 := {\displaystyle \frac {1}{1876463765}} ( - 26177 + 
 2292\,\sqrt{59})(2996295\,x^{2} - 256330\,y^{2} + 24928\,y\,x\,
  \sqrt{59} \\
  &\mbox{} + 14372\,y\,\sqrt{59} - 371968\,x\,\sqrt{59} + 306120\,x
  ^{2}\,\sqrt{59} + 956794 + 183424\,y\,x \\
  &\mbox{} + 77308\,\sqrt{59} + 346616\,y - 3822204\,x) 
  \end{eqnarray*}
  Notice that two conics $f_2(x,y)=0$ and $g_2(x,y)=0$ intersect at (1,1)
  with intersection multiplicity 3 to make $A_8$
with one more transverse intersection which makes an inner $A_2$.

  \vspace{.5cm}\noindent
  3. We consider $\Si_{4,1}:=[3A_5],\, \Si_{4,2}:=[2A_5+E_6],\, \Si_{4,3}:=[A_5+2E_6]$.
 For $[2A_5+E_6]$, we use Method 2. For $[3A_5],\, [2E_6+A_5]$, we can
 use either Method 1 or 2. For the configuration $[3E_6]$, we have
 already  constructed in \cite{Okadual}. 

 \vspace{.4cm}\noindent
 3-1. A sextic with configuration $[3A_5]$ is given by Method 2 as:
 \begin{multline*}
 g(x,y)=f_2(x,y)^3+g_2(x,y)^2 h_2(x,y)\\
 f_2 := 1 + 2\,x^{2} - {\displaystyle \frac {25}{1407}} \,
 I\,x^{2}\,\sqrt{1407} - y^{2} + x\,y + {\displaystyle \frac {25}{
 1407}} \,I\,x\,\sqrt{1407} - 2\,x\\
 g_2 := {\displaystyle \frac {5}{4}} \,y^{2} - 
 {\displaystyle \frac {5}{4}} \,x\,y + {\displaystyle \frac {125}{
 5628}} \,I\,x^{2}\,\sqrt{1407} - {\displaystyle \frac {1}{4}} \,x
 ^{2} - {\displaystyle \frac {125}{5628}} \,I\,x\,\sqrt{1407} + 
 {\displaystyle \frac {5}{2}} \,x - {\displaystyle \frac {5}{4}}\\
 h_2 :=  - {\displaystyle \frac {155}{16}} \,y^{2} + 
 {\displaystyle \frac {2025}{3752}} \,I\,x\,y\,\sqrt{1407} - 
 {\displaystyle \frac {11}{2}} \,x\,y - {\displaystyle \frac {2025
 }{3752}} \,I\,y\,\sqrt{1407} + {\displaystyle \frac {131449}{7504
 }} \,x^{2} \\
 \mbox{} + {\displaystyle \frac {275}{2814}} \,I\,x^{2}\,\sqrt{
 1407} - {\displaystyle \frac {275}{2814}} \,I\,x\,\sqrt{1407} - 
 {\displaystyle \frac {138953}{3752}} \,x + {\displaystyle \frac {
 138953}{7504}}  
 \end{multline*}
 Another example is obtained using Method 1:
 \begin{multline*}
 f := 27800 + 7600\,y^{5} - 15200\,y^{3} + 86200\,x^{3}\,y^{3} + 
 1650\,y^{4}\,x^{2} + 15050\,x^{6} - 3900\,y^{2}\,x^{4} \\
 \mbox{} - 86200\,x^{3}\,y + 58000\,y\,x^{4} + 58000\,y\,x^{2} + 
 59650\,x^{2} - 58600\,x + 7600\,y \\
 \mbox{} - 58600\,y^{4}\,x + 117200\,y^{2}\,x - 58000\,y^{3}\,x^{2
 } - 40600\,y\,x^{5} + 40600\,y^{4} - 62000\,y^{2} \\
 \mbox{} + 61900\,x^{4} - 6400\,y^{6} - 61300\,y^{2}\,x^{2} - 
 65200\,x^{3} - 6400\,y^{3}\,x + 3200\,y\,x \\
 \mbox{} + 65200\,y^{2}\,x^{3} + 3200\,y^{5}\,x - 40600\,x^{5} 
 \end{multline*}

\noindent
 3-2. A sextics with $[2A_5+E_6]$ is given by the following.
We have used Method 1 and $2A_5$  are $(0,-1),\, (0,1)$ with
respective tangent cones $y\mp 1=0$ and $E_6$ is at $(1,0)$.
 
\begin{multline*}
{\displaystyle \frac {1664357}{25425}} \,x + 2\,y + 
{\displaystyle \frac {203063}{25425}} \,y^{2} - {\displaystyle 
\frac {135871}{2825}} \,y\,x + y^{6} - 4\,y^{3} + {\displaystyle 
\frac {18148}{565}} \,y\,x^{5} - {\displaystyle \frac {34736}{
2825}} \,y\,x^{2} \\
\mbox{} - {\displaystyle \frac {73177}{25425}} \,x^{6} - 
{\displaystyle \frac {88819}{25425}}  + 2\,y^{5} - 
{\displaystyle \frac {5222126}{25425}} \,x^{2} - {\displaystyle 
\frac {135871}{2825}} \,y^{5}\,x + {\displaystyle \frac {1664357
}{25425}} \,y^{4}\,x \\
\mbox{} + {\displaystyle \frac {6641603}{25425}} \,x^{3} - 
{\displaystyle \frac {139669}{25425}} \,y^{4} - {\displaystyle 
\frac {3903338}{25425}} \,x^{4} - {\displaystyle \frac {4882462}{
25425}} \,y^{4}\,x^{2} + {\displaystyle \frac {39260}{1017}} \,x
^{5} \\
\mbox{} + {\displaystyle \frac {271742}{2825}} \,y^{3}\,x - 193\,
y^{3}\,x^{3} - {\displaystyle \frac {3328714}{25425}} \,y^{2}\,x
 - {\displaystyle \frac {471008}{2825}} \,y\,x^{4} + 193\,x^{3}\,
y \\
\mbox{} - {\displaystyle \frac {6641603}{25425}} \,y^{2}\,x^{3}
 + {\displaystyle \frac {1122732}{2825}} \,y^{2}\,x^{2} + 
{\displaystyle \frac {34736}{2825}} \,y^{3}\,x^{2} - 
{\displaystyle \frac {337334}{25425}} \,y^{2}\,x^{4}
=0
\end{multline*}

\noindent
 3-3. A sextic with $[2E_6+A_5]$. We put $2 E_6$ at $(0,\pm 1)$
with respective tangent cones $y=\pm 1$ and we put $A_5$ at $(1,0))$.
Take $g=f_2^3+g_2h_2$ where
 \begin{multline*}
 f_2 :=  - {\displaystyle \frac {667}{676}} \,y^{2} - 
 {\displaystyle \frac {135}{676}} \,I\,y^{2}\,\sqrt{3} - 
 {\displaystyle \frac {53}{26}} \,x\,y + {\displaystyle \frac {15
 }{26}} \,I\,y\,x\,\sqrt{3} + y - x^{2} + x,\quad
 g_2 := y^{2}\\
 h_2 := {\displaystyle \frac {409455}{228488}} \,x^{2} + 
 {\displaystyle \frac {915381}{742586}} \,x\,y + {\displaystyle 
 \frac {2175255}{1485172}} \,I\,y\,x\,\sqrt{3} + {\displaystyle 
 \frac {16605}{28561}} \,I\,\sqrt{3} - {\displaystyle \frac {25029
 }{57122}}  \\
 \mbox{} + {\displaystyle \frac {20915415}{38614472}} \,I\,y^{2}\,
 \sqrt{3} - {\displaystyle \frac {291195}{228488}} \,I\,x\,\sqrt{3
 } + {\displaystyle \frac {18225}{228488}} \,I\,x^{2}\,\sqrt{3} - 
 {\displaystyle \frac {262035}{228488}} \,I\,y\,\sqrt{3} \\
 \mbox{} - {\displaystyle \frac {373005}{228488}} \,x + 
 {\displaystyle \frac {282123}{228488}} \,y - {\displaystyle 
 \frac {15287373}{19307236}} \,y^{2} 
 \end{multline*}

 \vspace{.5cm}\noindent
 4. Sextics with $[2A_5+2A_2],\, [A_5+E_6+2A_2],\, [2E_6+2A_2]$.

 \vspace{.3cm}\noindent
 {\bf 4-1.} A sextic with $[2A_5+2A_2]$ with $2A_5$ as inner
 singularities. We use Method 2. We
 put one inner $A_2$ at $(1,0)$ and an outer $A_5$
at $O$. Take $g=f_2^3+g_2h_2$ where
\begin{eqnarray*}
&f_2 := 1 - {\displaystyle \frac {1}{3}} \,y\,x + 
{\displaystyle \frac {8203}{31752}} \,y^{2} + {\displaystyle 
\frac {1}{3}} \,y - {\displaystyle \frac {5444}{3969}} \,x + 
{\displaystyle \frac {1475}{3969}} \,x^{2}\\
&g_2 := 1 - y\,x + {\displaystyle \frac {21139}{31752}} \,
y^{2} + y - {\displaystyle \frac {14857}{7938}} \,x + 
{\displaystyle \frac {6919}{7938}} \,x^{2}\\
&h_2 :=  - 1 + y\,x - {\displaystyle \frac {24667}{31752}
} \,y^{2} + y + {\displaystyle \frac {1475}{3969}} \,x - 
{\displaystyle \frac {36115}{47628}} \,x^{2}
\end{eqnarray*}
\noindent
{\bf 4-2.}  A sextic with $[A_5+E_6+2A_2]$ with $E_6+2A_2$ as inner
 singularities. We put $E_6$ at $(1,0)$ as an inner singularity
and  an $A_5$ as an outer singularity at $O$. Take $g=f_2^3+g_2h_2$ where
\begin{eqnarray*}
&f_2 :=  - {\displaystyle \frac {1}{24}} \,y^{2} - y\,x + 
y + {\displaystyle \frac {13}{12}} \,x^{2} - {\displaystyle 
\frac {25}{12}} \,x + 1\\
&g_2 := {\displaystyle \frac {3}{8}} \,y^{2} - 
{\displaystyle \frac {3}{2}} \,y\,x + {\displaystyle \frac {3}{2}
} \,y + x^{2} - 2\,x + 1\\
&h_2 := {\displaystyle \frac {1}{8}} \,y^{2} - 
{\displaystyle \frac {1}{4}} \,y\,x - {\displaystyle \frac {37}{
32}} \,x^{2} + {\displaystyle \frac {9}{4}} \,x - 1
\end{eqnarray*}
Note that $g_2(x,y)=-1/16\, \left( 4\,x-3\,y+y\sqrt {3}-4 \right)  \left( -4\,x+4+3\,y+y
\sqrt {3} \right)$.

\noindent
{\bf 4-3.} A sextic with $[2E_6+2A_2]$ with $2E_6$ as inner
 singularities at $O$ and $(1,0)$.  $2\,A_2$ are at $(0,1),\,(1,1)$.
Take $g=f_2^3+g_2h_2$ where
\begin{multline*}
f_2 := {\displaystyle \frac {14}{9}} \,y^{2} + 
{\displaystyle \frac {2}{27}} \,I\,y^{2}\,\sqrt{42} + y\,x - 
{\displaystyle \frac {5}{9}} \,y - {\displaystyle \frac {2}{27}} 
\,I\,y\,\sqrt{42} - {\displaystyle \frac {1}{9}} \,x^{2} + 
{\displaystyle \frac {1}{9}} \,x,\quad
g_2 := y^{2}\\
h_2 :=  - {\displaystyle \frac {23}{12}} \,y^{2} - 
{\displaystyle \frac {5}{9}} \,I\,y^{2}\,\sqrt{42} - 7\,y\,x - 
{\displaystyle \frac {2}{3}} \,I\,y\,x\,\sqrt{42} + 
{\displaystyle \frac {1}{6}} \,y + {\displaystyle \frac {8}{9}} 
\,I\,y\,\sqrt{42} - {\displaystyle \frac {11}{3}} \,x^{2} + 
{\displaystyle \frac {11}{3}} \,x \\
+ {\displaystyle \frac {2}{3}} \,I\,x\,\sqrt{42}
\mbox{} + {\displaystyle \frac {3}{4}}  - {\displaystyle \frac {1
}{3}} \,I\,\sqrt{42} 
\end{multline*}

\vspace{.5cm}\noindent
5.  A sextic with $[A_5+4A_2],\, [E_6+4A_2]$.

\noindent
{\bf 5-1.} A sextic with $[A_5+4A_2]$ is given by Method 2. 
We put  $4\, A_2$ as inner singularities
and an $A_5$ at $O$ as an outer singularity. Take $g=f_2^3+g_2h_2$ where
\begin{multline*}
f_2(x,y) := 1 + {\displaystyle \frac {33863}{1650}} \,x + y
 - {\displaystyle \frac {33863}{1650}} \,x^{2} + {\displaystyle 
\frac {23404966}{226875}} \,x\,y - {\displaystyle \frac {12222527
}{907500}} \,y^{2}\\
g_2(x,y) := 2 + {\displaystyle \frac {13727}{550}} \,x + 7\,y
 - {\displaystyle \frac {14277}{550}} \,x^{2} + {\displaystyle 
\frac {43781281}{302500}} \,x\,y - {\displaystyle \frac {1217817
}{75625}} \,y^{2}\\
h_2(x,y) :=  - {\displaystyle \frac {1}{4}}  - 
{\displaystyle \frac {2517}{275}} \,x + y + {\displaystyle 
\frac {9243}{1100}} \,x^{2} - {\displaystyle \frac {1667401}{
151250}} \,x\,y + {\displaystyle \frac {209923}{151250}} \,y^{2}
\end{multline*}

\vspace{.4cm}\noindent
{\bf 5-2. } A sextic with $[E_6+4A_2] $ is obtained by Method 2.
We put $2\,A_2$ at $(0,1),\, (1,0)$ as outer singularities and an inner
$E_6$
at $(1,1)$. Take $g=f_2^3+g_2h_2$ where
\begin{multline*}
f_2 := ({\displaystyle \frac {163}{171}}  - 
{\displaystyle \frac {1}{171}} \,\sqrt{854})\,y^{2} + (( - 
{\displaystyle \frac {221}{152}}  - {\displaystyle \frac {1}{342}
} \,\sqrt{854})\,x - {\displaystyle \frac {683}{1368}}  + 
{\displaystyle \frac {1}{114}} \,\sqrt{854})\,y \\
\mbox{} + ({\displaystyle \frac {31}{152}}  + {\displaystyle 
\frac {1}{342}} \,\sqrt{854})\,x^{2} + x - {\displaystyle \frac {
31}{152}}  - {\displaystyle \frac {1}{342}} \,\sqrt{854}
\end{multline*}
\begin{multline*}
g_2 := y^{2} + ( - {\displaystyle \frac {7}{8}} \,x - 
{\displaystyle \frac {9}{8}} )\,y + {\displaystyle \frac {1}{8}} 
\,x^{2} + {\displaystyle \frac {5}{8}} \,x + {\displaystyle 
\frac {1}{4}}\\
h_2 := ( - {\displaystyle \frac {2179}{4332}}  + 
{\displaystyle \frac {4}{1083}} \,\sqrt{854})\,y^{2} + ((
{\displaystyle \frac {1350}{361}}  - {\displaystyle \frac {49}{
1083}} \,\sqrt{854})\,x - {\displaystyle \frac {2036}{1083}}  + 
{\displaystyle \frac {10}{361}} \,\sqrt{854})\,y \\
\mbox{} + ( - {\displaystyle \frac {53}{1083}} \,\sqrt{854} - 
{\displaystyle \frac {131}{1444}} )\,x^{2} + ({\displaystyle 
\frac {29}{361}} \,\sqrt{854} - {\displaystyle \frac {1655}{722}
} )\,x - {\displaystyle \frac {34}{1083}} \,\sqrt{854} + 
{\displaystyle \frac {1997}{1444}}  
\end{multline*}

\section{Further application}
The semi-torus method is also very useful to construct other plane
curves. We give an example.
We can construct a curve $C$ of degree 10 with 10 $A_4$
singularities  where $8$ of them are inner and two $A_4$ are outer and
we put them
at $(0,\pm 1)$.
\begin{eqnarray*}
&f(x,y)=f_2(x,y)^5\,+\,g_4(x,y)^2\, h_2(x,y)\\
&f_2(x,y)= - {\displaystyle \frac {9}{4}} \,y^{2} + {\displaystyle \frac {
19}{4}} \,x^{2} + x + {\displaystyle \frac {13}{4}} \\
&g_4(x,y) :=  - {\displaystyle \frac {105}{32}} \,y^{4} - 
{\displaystyle \frac {199}{16}} \,x^{4} - 8\,x^{3} + 
{\displaystyle \frac {165}{16}} \,y^{2} - {\displaystyle \frac {
21}{4}} \,x - {\displaystyle \frac {1}{24}} \,x^{3}\,\sqrt{30} - 
{\displaystyle \frac {661}{32}} \,x^{2} + {\displaystyle \frac {1
}{12}} \,x^{2}\,\sqrt{30} \\
&\mbox{} + {\displaystyle \frac {405}{32}} \,y^{2}\,x^{2} + 
{\displaystyle \frac {15}{4}} \,y^{2}\,x - {\displaystyle \frac {
1}{24}} \,x^{4}\,\sqrt{30} - {\displaystyle \frac {257}{32}} \\
&h_2(x,y)={\displaystyle \frac {15}{4}} \,y^{2} - {\displaystyle \frac {19
}{4}}  - {\displaystyle \frac {813}{16}} \,x^{2} - 2\,( - 
{\displaystyle \frac {661}{32}}  + {\displaystyle \frac {1}{12}} 
\,\sqrt{30})\,x^{2} - 2\,x
\end{eqnarray*}
Consider the normal form of $A_4$ singularity $P\in P$:
$v^2+u^5=0$.
The ideal of adjuction $\cI_{P,k,10}$ (see \cite{Artal, OkaAtlas} for 
definition) is defined as
\[
 \cI_{P,10,9}=\langle u^2,v\rangle,
\quad \cI_{P,k,10}=\langle u,v\rangle\,\, K=7,8
\]
and the explicit computation of
$\si_k: \cO(k-3)\to \oplus \cO_P/\cI_{P,k,10}$
shows that the Alexander polynomial of $C$ is trivial.
Thus together with a generic curve of $(5,2)$-torus $C'$,
 we have a Zariski pair $(C,C')$
which are distinguished by their Alexander polynomials.


\def\cprime{$'$} \def\cprime{$'$} \def\cprime{$'$} \def\cprime{$'$}
  \def\cprime{$'$}

\end{document}